\documentclass[a4paper,oneside,english]{amsart}
\usepackage{setspace}
\usepackage{amssymb}

\makeatletter

 \theoremstyle{remark}    
 \newtheorem*{acknowledgement*}{Acknowledgements} 
 \newcommand{\lyxaddress}[1]{
   \par {\raggedright #1 
   \vspace{1.4em}
   \noindent\par}
 }
 \theoremstyle{plain}    
 \newtheorem*{lem*}{Lemma} 
 \theoremstyle{remark}
 \newtheorem*{rem*}{Remark}
 \theoremstyle{plain}    
 \newtheorem{thm}{Theorem} 

\usepackage{babel}
\makeatother
\begin{document}

\title{Deformed commutators on quantum group module-algebras}

\author{A. O. Garc\'{\i}a}

\lyxaddress{Max-Planck-Institut f\"ur Physik, F\"ohringer Ring 6, 80805 M\"unchen,
Germany.}

\email{ariel@mppmu.mpg.de}



\begin{abstract}
We construct quantum commutators on module-algebras of quasi-tri\-an\-gu\-lar
Hopf algebras. These are quantum-group covariant, and have generalized
antisymmetry and Leibniz properties. If the Hopf algebra is triangular
they additionally satisfy a generalized Jacobi identity, turning the
module-algebra into a quantum-Lie algebra. 
\end{abstract}
\maketitle
The purpose of this short communication is to present a quantum commutator
structure which appears naturally on any module algebra $A$ of a
quantum group $H$. In section \ref{sec:q-commutator} we write down
the main properties we require from a generalized commutator on a
quantum group module-algebra, and we give its definition. In section
\ref{sec:commutator_properties} we prove a theorem collecting the
main properties of this algebraic structure. Finally, in section \ref{sub:qPlane_example}
we develop an example, showing some explicit calculations for the
reduced $SL_{q}(2,\mathbb{C})$ quantum plane. We refer the reader
to the Appendix for notation and some basic facts on quasi-triangular
Hopf-algebras.

\section{The $q$-commutator\label{sec:q-commutator}}

Let $H$ be a quasi-triangular Hopf algebra. Take $A$ some $H$-module-algebra
(a left one, say). As usual, we will denote the action of $h\in H$
on $a\in A$ by $h\triangleright a$, and the coproduct using the
Sweedler notation $\Delta h=h_{1}\otimes h_{2}$. Being a left-module-algebra,
of course $h\triangleright \left(ab\right)=\left(h_{1}\triangleright a\right)\left(h_{2}\triangleright b\right)$.
As our main goal is to define a covariant commutator for which some
generalized Leibniz rule holds on both variables, a natural way to
start is proposing a deformation of the usual $\left[a,b\right]=ab-ba$
structure valid on any associative algebra. The deformation we start
with is\begin{eqnarray}
\left[a\, ,\, b\right]_{\chi } & \equiv  & m\circ \left(1-\chi \right)(a\otimes b)\label{eq:q-commutator}\\
 & = & ab-m(\chi (a\otimes b))\nonumber 
\end{eqnarray}
Here $m$ is the product on $A$ and the linear map\[
\chi \, :\, A\otimes A\longmapsto A\otimes A\; ,\]
which replaces the standard transposition operator $\tau $, needs
to be determined. Later on, we will sometimes use the generic decomposition\begin{equation}
\chi (a\otimes b)=\sum _{i}\sigma _{a}^{i}(b)\otimes e_{i}\qquad \qquad \left\{ e_{i}\right\} \; \textrm{vector space basis of}\; V\label{eq:explicit_form_for_braiding}\end{equation}
Clearly, the maps $\sigma ^{i}$ have to be linear in both $a$ and
$b$.

The most basic property we require the commutators to satisfy is some
adequate generalization of the Leibniz rule, on both variables. Such
a rule means that commuting the first (say) variable $a$ to the right
through a product $bc$ must be equivalent to commuting it in two
steps, first through $b$ and then through $c$. Expressed in terms
of the map $\chi $, which generalizes and deforms the permutation,
this would read\begin{eqnarray*}
\chi \circ \left(1\otimes m\right) & = & \left(m\otimes 1\right)\circ \left(1\otimes \chi \right)\circ \left(\chi \otimes 1\right)\\
\chi \circ \left(m\otimes 1\right) & = & \left(1\otimes m\right)\circ \left(\chi \otimes 1\right)\circ \left(1\otimes \chi \right)
\end{eqnarray*}
where the second relation come from commuting the second variable
$c$ to the left through a product $ab$. 

Note now the analogy between the above conditions and the ones required
on the braiding \cite{Majid}\[
\chi _{V,W}\, :\, V\otimes W\longmapsto W\otimes V\]
 of a braided monoidal category. These are\begin{eqnarray*}
\chi _{V,W\otimes U} & = & \left(1\otimes \chi _{V,U}\right)\circ \left(\chi _{V,W}\otimes 1\right)\\
\chi _{V\otimes W,U} & = & \left(\chi _{V,U}\otimes 1\right)\circ \left(1\otimes \chi _{W,U}\right)
\end{eqnarray*}
and illustrate the fact that moving an element of $V$ to the right
through $W\otimes U$ (resp. an element of $U$ to the left through
$V\otimes W$) should produce the same result if it is done in one
or two steps. Note that $V$, $W$ and $U$ are not even vector spaces
in the general case, and that our map $\chi $ acts on an algebra.
However, remembering the standard result that shows that the category
of $H$-modules of a quasi-triangular Hopf-algebra $H$ is braided
(see \cite{Majid}, for instance), we take here the same braiding
as an Ansatz and we will show in the next section that it satisfies
the required conditions.

Concretely, we take 

\begin{eqnarray}
\chi (a\otimes b) & \equiv  & \left(R_{2}\triangleright b\right)\otimes \left(R_{1}\triangleright a\right)\label{eq:braiding_ansatz}
\end{eqnarray}
where\[
R\equiv R_{1}\otimes R_{2}\]
is the $R$-matrix of $H$ (c.f. Appendix). Of course, a generic sum
of the type $R=\sum _{k}R_{1}^{k}\otimes R_{2}^{k}$ is understood.
Note that we could also use the second quasi-triangular structure
$\bar{R}$, obtaining a map $\bar{\chi }$ which will differ from
$\chi $ unless $H$ is triangular. As it is easy to see from the
definition of $\bar{R}$, this second map $\bar{\chi }$ is the inverse
of the first one,\[
\bar{\chi }\circ \chi =\chi \circ \bar{\chi }=1\]
The properties of $R$ imply now\begin{eqnarray*}
\chi (a\otimes \mathbf{1}) & = & \mathbf{1}\otimes a\\
\chi (\mathbf{1}\otimes a) & = & a\otimes \mathbf{1}
\end{eqnarray*}
and therefore\[
[\mathbf{1}\, ,\, a]_{\chi }=[a\, ,\, \mathbf{1}]_{\chi }=0\qquad \forall a\in A\]
However, note that in general it will be\[
[a\, ,\, a]_{\chi }\neq 0\]
because $\chi (a\otimes a)=\left(R_{2}\triangleright a\right)\otimes \left(R_{1}\triangleright a\right)$
is a priori different from $a\otimes a$.

\section{Properties of the commutator\label{sec:commutator_properties}}

\subsection{$q$-Leibniz rules}

As was the aim when defining the deformed commutator, we have

\begin{lem*}
The map $\left[\: ,\, \right]_{\chi }$ has a Leibniz property on
the second variable reading

\begin{eqnarray}
\left[a\, ,\, bc\right]_{\chi } & = & \left[a\, ,\, b\right]_{\chi }\, c+\sigma _{a}^{i}(b)\left[e_{i}\, ,\, c\right]_{\chi }\label{eq:Leibniz_rule_2nd}
\end{eqnarray}
or, equivalently,\begin{eqnarray}
\chi \left(a\otimes bc\right) & = & \left(m\otimes 1\right)\left(1\otimes \chi \right)\left(\chi \left(a\otimes b\right)\otimes c\right)\; .\label{eq:Leibniz_rule_2nd_formal}
\end{eqnarray}
The corresponding equations for the Leibniz rule on the first variable
are\begin{eqnarray}
\left[ab\, ,\, c\right]_{\chi } & = & \left[a\, ,\, \sigma _{b}^{i}(c)\right]_{\chi }\, e_{i}+a\left[b\, ,\, c\right]_{\chi }\label{eq:Leibniz_rule_1st}
\end{eqnarray}
or, equivalently,\begin{eqnarray}
\chi \left(ab\otimes c\right) & = & \left(1\otimes m\right)\left(\chi \otimes 1\right)\left(a\otimes \chi \left(b\otimes c\right)\right)\; .\label{eq:Leiniz_rule_1st_formal}
\end{eqnarray}

\end{lem*}
The equivalency between, say, (\ref{eq:Leibniz_rule_2nd}) and (\ref{eq:Leibniz_rule_2nd_formal})
is straightforward keeping in mind that $\sigma _{a}^{i}(b)\otimes e=\chi (a\otimes b)$
and the definition (\ref{eq:q-commutator}). Using the explicit notation
(\ref{eq:explicit_form_for_braiding}), the above properties translate
into\begin{eqnarray*}
\sigma _{a}^{i}(bc) & = & \sigma _{a}^{i'}(b)\, \sigma _{e_{i'}}^{i}(c)
\end{eqnarray*}
and\begin{eqnarray*}
\sigma _{ab}^{j}(c)\otimes e_{j} & = & \sigma _{a}^{i}\left(\sigma _{b}^{i'}(c)\right)\otimes e_{i}e_{i'}
\end{eqnarray*}
respectively.

We only write down here the proof of (\ref{eq:Leibniz_rule_2nd_formal}),
the one of (\ref{eq:Leiniz_rule_1st_formal}) corresponds to a trivial
alteration of the former. Expand\begin{eqnarray*}
\chi \left(a\otimes bc\right) & = & \left(R_{2}\triangleright \left(bc\right)\right)\otimes \left(R_{1}\triangleright a\right)
\end{eqnarray*}
Considering (\ref{eq:Delta_R}), this gives\begin{eqnarray*}
\chi \left(a\otimes bc\right) & = & \left(m\otimes 1\right)\left(\Delta R_{2}\otimes R_{1}\right)\triangleright \left(b\otimes c\otimes a\right)\\
 & = & \tau \left(1\otimes m\right)\left(R_{13}R_{12}\triangleright \left(a\otimes b\otimes c\right)\right)
\end{eqnarray*}
Rewriting the action of $R_{12}$ in terms of $\chi $, and using
the trivial result $\tau \left(1\otimes m\right)=\left(m\otimes 1\right)\left(1\otimes \tau \right)\left(\tau \otimes 1\right)$,
we find\begin{eqnarray*}
\chi \left(a\otimes bc\right) & = & \left(m\otimes 1\right)\left(1\otimes \tau \right)\left(\tau \otimes 1\right)\left(R_{13}\triangleright \left(\tau \left[\chi \left(a\otimes b\right)\right]\otimes c\right)\right)\\
 & = & \left(m\otimes 1\right)\left(1\otimes \tau \right)\left(R_{23}\triangleright \left(\chi \left(a\otimes b\right)\otimes c\right)\right)\\
 & = & \left(m\otimes 1\right)\left(1\otimes \chi \right)\left(\chi \left(a\otimes b\right)\otimes c\right)
\end{eqnarray*}
which is the intended result.

\subsection{Covariance}

We will now prove

\begin{lem*}
The commutator $\left[\: ,\, \right]_{\chi }$ is quantum-group covariant,
in the sense that\begin{eqnarray}
h\triangleright [a,b]_{\chi } & = & \left[h_{1}\triangleright a\, ,\, h_{2}\triangleright b\right]_{\chi }\; .\label{eq:covariance}
\end{eqnarray}

\end{lem*}
Using the definition of the commutator and the quantum group action
properties,

\begin{eqnarray*}
h\triangleright [a,b]_{\chi } & = & h\triangleright \left(ab-\left(R_{2}\triangleright b\right)\left(R_{1}\triangleright a\right)\right)\\
 & = & \left(h_{1}\triangleright a\right)\left(h_{2}\triangleright b\right)-m\left[\left(\Delta h\, R^{\tau }\right)\triangleright \left(b\otimes a\right)\right]
\end{eqnarray*}
But according to (\ref{eq:R_intertwines_Delta}) we see that the last
term can be rewritten\begin{eqnarray*}
m\left[\left(\Delta h\, R^{\tau }\right)\triangleright \left(b\otimes a\right)\right] & = & m\, \tau \left[\left(\Delta ^{op}h\, R\right)\triangleright \left(a\otimes b\right)\right]\\
 & = & m\, \tau \left[\left(R\, \Delta h\right)\triangleright \left(a\otimes b\right)\right]\\
 & = & \chi \left[\Delta h\triangleright \left(a\otimes b\right)\right]
\end{eqnarray*}
Therefore\begin{eqnarray*}
h\triangleright [a,b]_{\chi } & = & m\, \left(1-\chi \right)\left[\Delta h\triangleright \left(a\otimes b\right)\right]\; ,
\end{eqnarray*}
which coincides with (\ref{eq:covariance}).

\subsection{$q$-Antisymmetry}

Generalizing the classical antisymmetry of a commutator, we now have

\begin{lem*}
The commutator $\left[\: ,\, \right]_{\chi }$ is $q$-antisymmetric,
this meaning

\begin{eqnarray}
[a,b]_{\chi } & = & -\left[\sigma _{a}^{i}(b)\, ,\, e_{i}\right]_{\bar{\chi }}\nonumber \\
 & = & -\left[\: ,\, \right]_{\bar{\chi }}\, \left(\chi (a\otimes b)\right)\label{eq:antisymmetry}
\end{eqnarray}

\end{lem*}
Note that in the RHS we have the deformed commutator $\left[\: ,\, \right]_{\bar{\chi }}$
given by the opposite quasi-triangular structure $\bar{R}$. The proof
is simply expressing the fact that $\bar{\chi }$ and $\chi $ are
inverse maps:\begin{eqnarray*}
[a,b]_{\chi } & = & m\, \left(1-\chi \right)\left(a\otimes b\right)\\
 & = & -m\, \left(1-\bar{\chi }\right)\, \chi \left(a\otimes b\right)
\end{eqnarray*}
If the quantum group $H$ is triangular, $\bar{R}=R$ and a same and
unique commutator appears in (\ref{eq:antisymmetry}).

\subsection{Conjugacy properties}

Let us now analyze the conjugacy properties of the commutator with
respect to a star operation on $A$. Assume\begin{eqnarray*}
\star _{H} & : & H\longmapsto H
\end{eqnarray*}
 is a Hopf-star on $H$, and \begin{eqnarray*}
\star _{A} & : & A\longmapsto A
\end{eqnarray*}
 is a compatible star \cite{CoGaTr-Stars} on $A$, in the sense that
\begin{equation}
h\triangleright \left(a^{\star _{A}}\right)=\left[\left(Sh\right)^{\star _{H}}\triangleright a\right]^{\star _{A}}\; .\label{eq:compatible_Hopf_stars}\end{equation}
Then we can analyze the conjugacy properties of the commutator. From
now on we drop the indexes on $\star $, as there is no confusion
possible.

\begin{lem*}
If $R$ is anti-real \cite{Majid}, meaning \begin{eqnarray}
R^{\star } & = & R^{-1}\; ,\label{eq:antireal_R}
\end{eqnarray}
then\begin{eqnarray*}
[a,b]_{\chi }^{\star } & = & \left[b^{\star }\, ,\, a^{\star }\right]_{\bar{\chi }}\; .
\end{eqnarray*}
 For a real $R$, i.e. such that\begin{eqnarray}
R^{\star } & = & \tau \left(R\right)\; ,\label{eq:real_R}
\end{eqnarray}
the result is\begin{eqnarray*}
[a,b]_{\chi }^{\star } & = & \left[b^{\star }\, ,\, a^{\star }\right]_{\chi }\; .
\end{eqnarray*}

\end{lem*}
The quantum plane example shown in section \ref{sub:qPlane_example}
corresponds to the first possibility. The proof goes as follows:\begin{eqnarray*}
[a,b]_{\chi }^{\star } & = & b^{\star }a^{\star }-\left(R_{1}\triangleright a\right)^{\star }\left(R_{2}\triangleright b\right)^{\star }\; .
\end{eqnarray*}
Considering first (\ref{eq:compatible_Hopf_stars}), and using next
that $\left(S\otimes S\right)R=R$, we obtain\begin{eqnarray*}
[a,b]_{\chi }^{\star } & = & b^{\star }a^{\star }-\left(\left(SR_{1}\right)^{\star }\triangleright a^{\star }\right)\left(\left(SR_{2}\right)^{\star }\triangleright b^{\star }\right)\\
 & = & b^{\star }a^{\star }-m\left[R^{\star }\triangleright \left(a^{\star }\otimes b^{\star }\right)\right]\\
 & = & m\left[1-\tau \circ \left(\tau \left(R^{\star }\right)\triangleright \, \cdot \, \right)\right]\left(b^{\star }\otimes a^{\star }\right)
\end{eqnarray*}
For a real $R$ (resp. anti-real), $\tau \left(R^{\star }\right)=R$
(resp. $=\bar{R}$) and the lemma follows.

\subsection{Quantum Lie algebra structure and Jacobi identities}

Having defined a generalized commutator with Leibniz and antisymmetry
properties, we could now inquire about the relationship between this
structure and the one provided by a quantum Lie algebra \cite{BuIsOg}.
Following this reference, a quantum Lie algebra is defined by relations\begin{equation}
e_{i}e_{j}-\sigma _{ij}^{mk}\, e_{m}e_{k}=C_{ij}^{k}\, e_{k}\label{eq:qLie_commutator}\end{equation}
among vector space generators $\left\{ e_{i}\right\} $ of the space.
The matrix $\sigma _{ij}^{mk}$ should satisfy a Yang-Baxter equation,
and the structure constants $C_{ij}^{k}$ have to obey equations (2),
(3), and (4) of \cite{BuIsOg}, corresponding to generalized Jacobi
and Leibniz properties. Comparing (\ref{eq:qLie_commutator}) with
(\ref{eq:q-commutator}) we see that we must take\[
\sigma _{ij}^{mk}\, e_{m}\otimes e_{k}=\chi (e_{i}\otimes e_{j})\]
and\[
C_{ij}^{k}\, e_{k}=\left[e_{i}\, ,\, e_{j}\right]_{\chi }\; .\]
Remark also that the Yang-Baxter equation (\ref{eq:YB}) implies for
$\chi $ the following relation:\begin{eqnarray}
\left(1\otimes \chi \right)\left(\chi \otimes 1\right)\left(1\otimes \chi \right) & = & \left(\chi \otimes 1\right)\left(1\otimes \chi \right)\left(\chi \otimes 1\right)\; .\label{eq:YB_for_braiding}
\end{eqnarray}
The proof is straightforward. Now using our (\ref{eq:YB_for_braiding}),
(\ref{eq:Leibniz_rule_2nd_formal}), and (\ref{eq:Leiniz_rule_1st_formal})
it is straightforward algebra to see that the conditions (3) and (4)
of \cite{BuIsOg} are satisfied.

Condition (2) of \cite{BuIsOg} corresponds to the Jacobi identity
of the quantum Lie algebra, and we have not yet analyzed such a property
for the commutators $\left[\: ,\, \right]_{\chi }$.

The usual Jacobi identity can, a priori, be generalized in several
possible ways. However, in order to maintain the parallel with the
$q$-Lie algebras of \cite{BuIsOg}, we take here the generalization

\begin{eqnarray}
\left[\left[\, \cdot \, ,\, \cdot \, \right]_{\chi },\, \cdot \, \right]_{\chi } & = & \left[\, \cdot \, ,\left[\, \cdot \, ,\, \cdot \, \right]_{\chi }\right]_{\chi }+\left[\left[\, \cdot \, ,\, \cdot \, \right]_{\chi },\, \cdot \, \right]_{\chi }\circ \left(1\otimes \chi \right)\label{eq:q-Jacobi}
\end{eqnarray}
which corresponds to their equation (2). After using the Leibniz properties,
(\ref{eq:q-Jacobi}) translates into\begin{eqnarray*}
\left\{ 1-\left(\chi \otimes 1\right)\left(1\otimes \chi \right)\right\} \left\{ 1-\left(\chi \otimes 1\right)\right\}  & = & \left\{ 1-\left(1\otimes \chi \right)\left(\chi \otimes 1\right)\right\} \left\{ 1-\left(1\otimes \chi \right)\right\} \\
 &  & +\left\{ 1-\left(\chi \otimes 1\right)\left(1\otimes \chi \right)\right\} \left\{ 1-\left(\chi \otimes 1\right)\right\} \left(1\otimes \chi \right)
\end{eqnarray*}
Making use of the Yang-Baxter equation for $\chi $ (\ref{eq:YB_for_braiding})
we get:\begin{eqnarray*}
0 & = & \left\{ \left(\chi \otimes 1\right)-\left(1\otimes \chi \right)\left(\chi \otimes 1\right)\right\} \left(1-\left(1\otimes \chi ^{2}\right)\right)
\end{eqnarray*}
 Therefore, the Jacobi identity is satisfied only in the case $\chi ^{2}=1$,
i.e., if $H$ is a triangular Hopf algebra.\\

All the above results can be collected in a

\begin{thm}
\label{thm:theorem}Let $A$ be a left-module-algebra of a quasi-triangular
Hopf algebra $H$, and take $\left[\: ,\, \right]_{\chi }$ the quantum
commutator on $A$ defined by (\ref{eq:q-commutator}) and (\ref{eq:braiding_ansatz}).
Then $\left[\: ,\, \right]_{\chi }$ is quantum-group covariant and
has generalized antisymmetry and Leibniz properties. If the $H$ is
triangular, then they additionally satisfy a generalized Jacobi identity,
turning the module-algebra into a quantum-Lie algebra.
\end{thm}

\section{The quantum plane example\label{sub:qPlane_example}}

Take the quantum plane algebra $A$ generated by $x$ and $y$ such
that

\[
xy=q\, yx\; ,\qquad \qquad q\in \mathbb{C}\: ,\quad q\ne 0\; .\]
On $A$ we have the action \cite{CoGaTr-qPlaneE} of the quantum enveloping
algebra $H=U_{q}(sl(2,\mathbb{C}))$ generated by $K$, $K^{-1}$,
$X_{+}$, $X_{-}$with relations\begin{eqnarray*}
KX_{\pm } & = & q^{\pm 2}\, X_{\pm }K\\
\left[X_{+}\, ,\, X_{-}\right] & = & \frac{1}{\left(q-q^{-1}\right)}\left(K-K^{-1}\right)
\end{eqnarray*}
Additionally one can take the complex parameter $q$ to be a root
of unit, $q^{N}=1$ for some (odd) integer $N$. In such a case one
can get non-trivial finite dimensional algebras by taking the quotient
of the above ones by the following ideals:\begin{eqnarray*}
x^{N} & = & \mathbf{1}\\
y^{N} & = & \mathbf{1}
\end{eqnarray*}
and\begin{eqnarray*}
K^{N} & = & \mathbf{1}\\
X_{\pm }^{N} & = & 0
\end{eqnarray*}
Of course now $K^{-1}=K^{N-1}$. To be concrete, we take the value
$N=3$, thus $q^{3}=1$. In this case the $R$-matrix of $U_{q}(sl(2,\mathbb{C}))$
is given by \cite{CoGaTr-qPlaneE}\begin{equation}
R=\frac{1}{3}\, R_{K}R_{X}\; ,\label{eq:qPlane_R_matrix}\end{equation}
where\begin{eqnarray*}
R_{K} & = & \mathbf{1}\otimes \mathbf{1}+\left(\mathbf{1}\otimes K+K\otimes \mathbf{1}\right)+\left(\mathbf{1}\otimes K^{2}+K^{2}\otimes \mathbf{1}\right)\\
 &  & +\: q^{2}\left(K\otimes K^{2}+K^{2}\otimes K\right)+q\, K\otimes K+q\, K^{2}\otimes K^{2}\\
R_{X} & = & \mathbf{1}\otimes \mathbf{1}+\left(q-q^{-1}\right)\, X_{-}\otimes X_{+}+3q\, X_{-}^{2}\otimes X_{+}^{2}
\end{eqnarray*}
Applying formula (\ref{eq:braiding_ansatz}) we can calculate the
following elementary $\chi $'s:\begin{eqnarray*}
\chi (x\otimes x) & = & q^{2}\, x\otimes x\\
\chi (y\otimes x) & = & q\, x\otimes y\\
\chi (x\otimes y) & = & q\, y\otimes x+\left(q^{2}-1\right)x\otimes y\\
\chi (y\otimes y) & = & q^{2}\, y\otimes y
\end{eqnarray*}
The quantum plane algebra can be extended in a covariant way introducing
derivative operators $\partial _{x}$ and $\partial _{y}$ \cite{WessZumino}.
We refer the reader to \cite{CoGaTr-qPlaneE} for the complete algebraic
structure. Including these derivatives, the braiding is:\begin{eqnarray*}
\chi (\partial _{x}\otimes x) & = & q\, x\otimes \partial _{x}\\
\chi (\partial _{y}\otimes x) & = & q^{2}\, x\otimes \partial _{y}\\
\chi (\partial _{x}\otimes y) & = & q^{2}\, y\otimes \partial _{x}\\
\chi (\partial _{y}\otimes y) & = & (q-1)x\otimes \partial _{x}+q\, y\otimes \partial _{y}
\end{eqnarray*}
\begin{eqnarray*}
\chi (x\otimes \partial _{x}) & = & q\, \partial _{x}\otimes x+(q^{2}-q)\partial _{y}\otimes y\\
\chi (y\otimes \partial _{x}) & = & q^{2}\, \partial _{x}\otimes y\\
\chi (\partial _{x}\otimes \partial _{x}) & = & q^{2}\, \partial _{x}\otimes \partial _{x}\\
\chi (\partial _{y}\otimes \partial _{x}) & = & (q^{2}-1)\partial _{y}\otimes \partial _{x}+q\, \partial _{x}\otimes \partial _{y}
\end{eqnarray*}
\begin{eqnarray*}
\chi (x\otimes \partial _{y}) & = & q^{2}\, \partial _{y}\otimes x\qquad \qquad \qquad \quad \\
\chi (y\otimes \partial _{y}) & = & q\, \partial _{y}\otimes y\\
\chi (\partial _{x}\otimes \partial _{y}) & = & q\, \partial _{y}\otimes \partial _{x}\\
\chi (\partial _{y}\otimes \partial _{y}) & = & q^{2}\, \partial _{y}\otimes \partial _{y}
\end{eqnarray*}
Using the relations between derivatives and coordinates found in \cite{CoGaTr-qPlaneE},
we can now display a few non-trivial commutators. We have, for instance,
\begin{eqnarray*}
\left[x,x\right]_{\chi } & = & x^{2}-m\left(\chi \left(x\otimes x\right)\right)\\
 & = & \left(1-q^{2}\right)x^{2}\\
 &  & \\
\left[\partial _{x},x\right]_{\chi } & = & \partial _{x}x-m\left(\chi \left(\partial _{x}\otimes x\right)\right)\\
 & = & \mathbf{1}+\left(q^{2}-q\right)x\partial _{x}+\left(q^{2}-1\right)y\partial _{y}\\
 &  & \\
\left[x,\partial _{x}\right]_{\chi } & = & x\partial _{x}-m\left(\chi \left(x\otimes \partial _{x}\right)\right)\\
 & = & -q^{2}\mathbf{1}
\end{eqnarray*}

\begin{rem*}
Note that one could think about using the matrix representation of
the reduced quantum plane at $q^{3}=1$ as a way to define commutators.
Taking the explicit $3\times 3$ matrices \cite{CoGaTr-qPlaneL,CoGaTr-qPlaneE}

\[
\mathbf{x}=\left(\begin{array}{ccc}
 1 & 0 & 0\\
 0 & q^{-1} & 0\\
 0 & 0 & q^{-2}\end{array}
\right)\qquad \qquad \mathbf{y}=\left(\begin{array}{ccc}
 0 & 1 & 0\\
 0 & 0 & 1\\
 1 & 0 & 0\end{array}
\right)\]
we see that our above deformed commutator has nothing to do with the
commutator of these matrices. In fact $\left[\mathbf{x},\mathbf{x}\right]=0$
(as matrices), whereas $\left[x,x\right]_{\chi }=\left(1-q^{2}\right)x^{2}$,
as we saw above. Of course, the point is that the commutator defined
using these matrices doesn't have the covariance property of our deformed
commutator.
\end{rem*}

\section{Concluding remarks}

The main results of this communication are collected in Theorem \ref{thm:theorem},
involving the existence of a covariant commutator structure on any
module-algebra of a quasi-triangular Hopf algebra. This commutator
turns the module-algebra into a quantum Lie algebra in the case that
the quantum group acting on it is triangular. The fact that the deformed
Jacobi identity (\ref{eq:q-Jacobi}) is obeyed only for a \emph{triangular}
Hopf algebra seems to be independent of the way we choose to generalize
the Jacobi identity.

\begin{acknowledgement*}
The author is deeply indebted to O. Ogievetsky and R. Coquereaux for
their comments and discussions, and gratefully acknowledges the Max-Planck-Ge\-sell\-schaft
for financial support.
\end{acknowledgement*}
\appendix

\section*{Appendix: quasi-triangular Hopf algebras}

We remember here that a \emph{quasi-triangular} Hopf algebra $H$
\cite{Majid} has, by definition, an element $R\in H\otimes H$ with
the following properties:

\begin{eqnarray}
\Delta ^{op}h & = & R\, \Delta h\, R^{-1}\label{eq:R_intertwines_Delta}
\end{eqnarray}
\begin{eqnarray}
(\Delta \otimes 1)R & = & R_{13}R_{23}\label{eq:Delta_R}\\
(1\otimes \Delta )R & = & R_{13}R_{12}\nonumber 
\end{eqnarray}
It follows that $R$ satisfies the Yang-Baxter equation\begin{eqnarray}
R_{12}R_{13}R_{23} & = & R_{23}R_{13}R_{12}\; .\label{eq:YB}
\end{eqnarray}
The algebra $H$ automatically has a second quasi-triangular structure
given by the related element\begin{equation}
\bar{R}=\tau \left(R^{-1}\right)\; ,\label{eq:R_bar_structure}\end{equation}
where $\tau $ is the permutation of tensor product factors. If both
$R$ and $\bar{R}$ coincide one says that the Hopf algebra $H$ is
in fact \emph{triangular}. Two additional basic properties of the
$R$-matrix which we need in our proofs are\begin{eqnarray*}
\left(\epsilon \otimes 1\right)R\: =\: \left(1\otimes \epsilon \right)R & = & \mathbf{1}\; ,
\end{eqnarray*}
\begin{eqnarray*}
\left(S\otimes S\right)R & = & R\; .
\end{eqnarray*}

\end{document}